\documentclass[letter,10pt]{article}

\usepackage{amsfonts}
\usepackage{amssymb}
\usepackage{amsthm}
\usepackage{amsmath}
\usepackage{latexsym}
\usepackage[all]{xy}
\usepackage{mathrsfs}
\usepackage{graphicx}

\theoremstyle{plain}\newtheorem{thm}{Theorem}[section]
\theoremstyle{plain}\newtheorem{lem}[thm]{Lemma}
\theoremstyle{plain}\newtheorem{prop}[thm]{Proposition}
\theoremstyle{plain}
\theoremstyle{definition}\newtheorem{defn}[thm]{Definition}
\theoremstyle{remark}\newtheorem{rem}[thm]{REMARK}
\theoremstyle{remark}\newtheorem{exa}[thm]{EXAMPLE}

\title{Reversible linear differential equations}
\author{Camilo Sanabria Malag\'on\footnote{partially supported by NSF grant CCF 0901175 and CCF 0952591.}}
\date{}

\begin{document}

\maketitle

\begin{abstract}
Let $\nabla$ be a meromorphic connection on a vector bundle over a compact Riemann surface $\Gamma$. An automorphism
$\sigma:\Gamma\rightarrow\Gamma$ is called a symmetry of $\nabla$ if the pullback bundle and the pullback connection can be identified
with $\nabla$. We study the symmetries by means of what we call the Fano Group; and, under the hypothesis that $\nabla$ has a unimodular
reductive Galois group, we relate the differential Galois group, the Fano group and the symmetries by means of an exact sequence.
\end{abstract}

\begin{quote}
\emph{``the solution of an intellectual problem comes about in a way not very different from what
happens when a dog carrying a stick in its mouth tries to get through a narrow
door: it will go on turning its head left and right until the stick slips through''} \hfill-Robert Musil
\end{quote}

\section{Introduction}

In 1900, G. Fano addressed the following problem~\cite{fano}: what are the consequences of algebraic relations between the solutions of a linear differential equation? The problem was apparently proposed to him by F. Klein. A particular concern was whether or not a linear differential equation with solutions satisfying a homogeneous polynomial can be ``solved in terms of linear equations of lower order''. This has been successfully studied by M. Singer, cf.~\cite{singer88}, and more recently by K.A. Nguyen, cf.~\cite{nguyen}.

Fano considered the group of projective automorphisms of the projective variety having the solutions of the differential equation as coordinate functions. This could be viewed as a primitive version of the differential Galois Group. Here we replace the Fano group with one slightly smaller: the group of projective automorphisms of the projective variety having as coordinate ring the $\mathbb{C}$-algebra generated by the solutions of the differential equation together with the $i$'th derivatives, $i\in\{1,\ldots,n-1\}$ (where $n$ is the order of the equation).

We treat the problem in terms of connections. Suppose we are given a ramified covering map $\phi: X'\longrightarrow  X$ of compact Riemann surfaces, together with a meromorphic vector bundle $E$ with a connection $\nabla$ over $ X$. We can use $\phi$ to pull back the bundle and the connection, obtaining $E'=\phi^*E$ together with $\nabla'=\phi^*\nabla$. In this context one obtains~\cite[Appendix B]{morales}  a natural injection
\[
\textrm{Gal}(\nabla')\longrightarrow \textrm{Gal}(\nabla),
\]
where the induced map on the Lie algebras is an isomorphism (i.e. the connected components of $1$ are isomorphic). Note that, a covering transformation $\sigma\in \textrm{Aut}_{\mathbb{C}}( X')$ of $\phi: X'\longrightarrow  X$ lifts to a horizontal automorphism of the vector bundle $(E',\nabla')$. The automorphisms lifting to the connection are called symmetries (of the connection). A concise exposition on symmetries is given in ~\cite{churchill}.

Conversely, suppose we begin with a connection on a meromorphic vector bundle over a Riemann surface which admits a symmetry (such equations are often called reversible). We can then consider the quotient Riemann surface with the canonically induced vector bundle and connection. This new connection has Galois group with bigger monodromy subgroup, but with identity component isomorphic to that of the original Galois group.

Checking for symmetries of a given connection over the Riemann sphere is quite easy (just consider permutations of singular points). Methods for revealing them in arbitrary contexts are far from simple. For example, one can consider the work by Dwork and Baldassarri~\cite{baldassarri}, ~\cite{baldassarri2}. The study of symmetries is intimately linked to the study of descent conditions and to the identification of pullbacks. These pullbacks, on their turn, are important in the classification of equations and in algorithmic implementations.

The purpose of this monograph is to suggest how the outer-automorphisms in the Fano group of the Galois group correspond to symmetries of the connection, and to give a proof in the case when the Galois group is reductive and unimodular and the connection is standard (Theorem \ref{theorem}).

\section{Algebraic justification}

\begin{rem}
The argument behind the proof of Theorem \ref{theorem} is of a geometric nature. Nevertheless, it is possible to describe the phenomena studied in this article algebraically. We do so informally in this section. For instance, the examples in the last section correspond more to the algebraic point of view than to the geometric one.
\end{rem}

We consider a field $k$ over $\mathbb{C}$ of transcendence degree one, together with a non-trivial derivation $v$. Recall that a non-trivial derivation is a non-zero additive map $v:k\rightarrow k$  satisfying the Leibnitz rule:
\[
v(fg)=v(f)g+fv(g).
\]
The collection of linear differential operators over $(k,v)$ forms a noncommutative ring extension $k[v]$ of $k$ with multiplication given by
\[
v\cdot f= v(f)+f\cdot v.
\]
As a consequence we see that every linear differential operator $L\in k[v]$ can be written in the form
\[
L = a_nv^n+a_{n-1}v^{n-1}+\ldots,+a_1v+a_0
\]
with $a_i\in k$ for $i\in\{0,1,\ldots,n\}$.

Let $\sigma^*\in \textrm{Aut}_{\mathbb{C}}(k)$ be an automorphism $k$ over $\mathbb{C}$. We can naturally lift $\sigma^*$ to an automorphism of $k[v]$ by setting
\[
\sigma^*v: f\mapsto v(\sigma^*f)
\]
Note that $\sigma^*v$ is a non-trivial derivation. In particular, since the $k$-module of derivations of $k$ has rank one over $k$, there is a non-zero element $v(\sigma)\in k$ such that:
\[
\sigma^* v=\frac{1}{v(\sigma)}v.
\]
One can define $\sigma^*$ to be a symmetry of $L$ if $\sigma^*L=f_\sigma\cdot L$ for some $f_\sigma\in k$. This means that the solutions to the equation $L(y)=0$ and to
$\sigma^*L(y)=0$ coincide; classically we say that $\sigma^*$ is a period of $L$.

Assume now that $v$ is invariant under $\sigma^*$ (i.e. that $v(\sigma)=1$). Then $v$ determines a derivation on the subfield $k^{\sigma^*}\subseteq k$ fixed by $\sigma^*$. In this case, if $\sigma^*$ is a symmetry of $L$, then $L$ has coefficients in $k^{\sigma^*}$ and $L$ restricts to a differential operator on $(k^{\sigma^*},v)$. Theorem \ref{theorem} says that the symmetries of $L$ manifest themselves as outer-automorphisms of the Galois group of $L$ over $(k,v)$.

\begin{rem}
In broad terms the geometric idea behind proof of the main result is the following. Let $ X$ be the projective algebraic curve with $\mathbb{C}( X)=k$. We fix a point $p\in X$. The differential module $k[v]/L$ defines a vector bundle with connection over $ X$~\cite[Chapter 2]{SvdP}. The Galois group $G$ of $L$ over $k$ acts on a fiber over $p$ of this vector bundle. Given $\sigma^*\in\textrm{Aut}_\mathbb{C}( X)$, it defines an automorphism $\sigma$ of $ X$. If $\sigma^*$ is a symmetry of $L$ we have two ways of identifying the fiber over $p$ and the fiber over $\sigma(p)$: by analytic extension through a path from $\sigma(p)$ to $p$, or via the symmetry $\sigma$. Jumping from one identification to the another amounts to acting on $G$ by an outer-automorphisms provided $L$ is standard (cf. Definition \ref{def0}). The geometric interpretation of $L$ having coefficients in $k^{\sigma^*}$ when $v$ is $\sigma^*$-invariant is: the connection defined by $k[v]/L$ descents to a connection over $ X^\sigma$. Explicitly it descends to the connection defined by $k^{\sigma^*}[v]/L$.
\end{rem}

\begin{rem}
Let us illustrate the standard hypothesis with an example from number fields. We take as base field $\mathbb{Q}(i)$. Consider the Galois extension of the polynomial 
\[
X^4-i=0;
\]
namely $\mathbb{Q}(i)(e^{\frac{\pi}{8}i})$. The Galois group is $\mathbb{Z}/4\mathbb{Z}$, and is generated by 
\[
e^{\frac{\pi}{8}i}\mapsto ie^{\frac{\pi}{8}i}.
\]
The arithmetic interpretation of the standard hypothesis is that
\[
\mathbb{Q}(i)(e^{\frac{\pi}{8}i})=\mathbb{Q}(e^{\frac{\pi}{8}i})
\]
which implies $\mathbb{Q}(i)\subseteq\mathbb{Q}(e^{\frac{\pi}{8}i})$. We can act on the extension $\mathbb{Q}(i)\subset\mathbb{Q}(i)(e^{\frac{\pi}{8}i})$ by complex conjugation. Complex conjugation transforms the polynomial $X^4-i$ into $X^4+i$. Notice that the Galois extension over $\mathbb{Q}(i)$  of these two polynomials is the same. As above the Galois group of $X^4+i$ is generated by $e^{\frac{-\pi}{8}i}\mapsto ie^{\frac{-\pi}{8}i}$, i.e. by
\[
e^{\frac{\pi}{8}i}\mapsto -ie^{\frac{\pi}{8}i}=i^3e^{\frac{\pi}{8}i}.
\]
We put everything together in the tower
\[
\mathbb{Q}\subset\mathbb{Q}(i)\subset\mathbb{Q}(i)(e^{\frac{\pi}{8}i})
\]
corresponding to the polynomial $X^8+1=(X^4-i)(X^4+i)$. The Galois group of the tower is the dihedral group of order $8$, and complex conjugation corresponds to the outer-automorphism of the cyclic subgroup inside the dihedral group. The discussion is summarized by the exact sequence
\[
1\longrightarrow \mathbb{Z}/4\mathbb{Z}\longrightarrow D_{2\cdot4}\longrightarrow\langle i\mapsto-i\rangle\longrightarrow 1.
\]
\end{rem}

\section{Setting and definitions}\label{secdef}

\begin{rem}\label{einsconv}
We remind that we use the Einstein summation convention for indices.
\end{rem}

Let $ X$ be a (connected) compact Riemann surface with field of meromorphic functions $k$.  Let
\[
\Pi:E\longrightarrow  X
\]
be an $n$-dimensional meromorphic vector bundle, with a meromorphic connection
\[
\nabla:\mathscr{E}\longrightarrow \Omega^1_\mathscr{M}\otimes_k \mathscr{E}
\]
where $\Omega^1_\mathscr{M}$ denotes the meromorphic differential forms over $ X$ and $\mathscr{E}$ the meromorphic sections of $\Pi$. We also denote by $\mathscr{T}X$ the vector fields of meromorphic tangent vectors to $ X$.
There is a natural map
\begin{eqnarray*}
\mathscr{T}X\otimes_k\Omega^1_\mathscr{M} & \longrightarrow & k\\
v\otimes\eta & \longmapsto & \eta(v)
\end{eqnarray*}
which canonically extends to
\begin{eqnarray*}
\mathscr{T}X\otimes_k\Omega^1_\mathscr{M}\otimes_k \mathscr{E} & \longrightarrow & \mathscr{E}\\
v\otimes\eta\otimes X &\longmapsto & \langle v, \eta\otimes X\rangle:= \eta(v)X.
\end{eqnarray*}
Given a meromorphic tangent vector field $v\in\mathscr{T}X$ we denote by $\nabla_v$  the derivation on $\mathscr{E}$
\[
\nabla_v(X)=\langle v,\nabla X \rangle.
\]

\begin{defn} Let $\sigma\in \textrm{Aut}_\mathbb{C}( X)$ (or equivalently $\sigma^*\in \textrm{Aut}_\mathbb{C}(k)$). We say that $\sigma$ is a \emph{symmetry} of $\nabla$ if there is a horizontal vector bundle morphism
\[
\tilde{\sigma}:(E,\nabla)\to(E,\nabla)
\]
lifting $\sigma$, i.e. $\tilde{\sigma}(\nabla_vX)=\nabla_{\sigma_*v}\tilde{\sigma}(X)$ and $\sigma\circ\Pi=\Pi\circ \tilde{\sigma}$.
\[
\xymatrix{
E\ar@{-->}[rr]^{\tilde{\sigma}}\ar[d]_\Pi & & E\ar[d]^\Pi \\
 X\ar[rr]_\sigma & &  X.
}
\]
The group of symmetries of $\nabla$ will  be denoted by $\textrm{Aut}_\nabla( X)$.
\end{defn}

\begin{rem}
A symmetry of $\nabla$ permutes the singular points.
\end{rem}

Let $(U,z)$ be a holomorphic chart of $ X$ centered at $p\in U$, where $U\subseteq  X$ is an open ball avoiding the singularities  of $\nabla$. $E$ is holomorphic and trivial above $U$ and $\nabla$ induces a holomorphic connection $\nabla'$ on $\Pi^{-1}U\rightarrow U$. Let $(\Pi^{-1}U,z,y^1,\ldots,y^n)$ be a trivializing chart of $E$. There exists a holomorphic horizontal frame $V_1,\ldots, V_n$ over $U$, i.e.
\[
v^i_j(z)=y^i(V_j(z))
\]
with $v^i_j(z)$ holomorphic in $U$ such that $\nabla'_{\frac{\partial}{\partial z}}V_j=0$ and $\det(v_i^j)(z)$ does not vanish in $U$  (see ~\cite{ince}).

\begin{defn}
The Fano group $G_F$ of $\nabla$ is the subgroup of $GL_n(\mathbb{C})$ fixing the homogeneous ideal in $\mathbb{C}[X^i_j]$ generated by the $G$-invariant homogeneous polynomials $P[X^i_j]\in\mathbb{C}[X^i_j]$ vanishing at $v^i_j$.
\end{defn}

\begin{rem}
The Fano group defined here differs from that considered in ~\cite{nguyen}. In fact our group $G_F$ contains the group $G^+$ used in ~\cite{nguyen}.  For instance, if the connection is standard (cf. Definition \ref{def0}) and the Galois group is finite, $G^+$ coincides with $G$ but $G_F$ may be larger (cf. Example \ref{ex1}). In ~\cite{fano}, G. Fano worked with automorphisms of projective varieties, so our approach is in the same spirit.
\end{rem}

\begin{rem}
It follows directly from the definition that there is a canonical inclusion of the representation in $GL_n(\mathbb{C})$ of the Galois group $G$ given by $v_j^i$ into the Fano group $G_F$. Let us make this remark more explicit.
\end{rem}

Fix $v\in\mathscr{T}X$, $v\not=0$, and a global meromorphic frame $(e_1,\ldots,e_n)$ of $\mathscr{E}$, i.e. an $n$-tuple of meromorphic sections such that on some Zariski open subset $ X'\subseteq X$, the $n$-tuple $(e_1(q),\ldots,e_n(q))$ is a basis of $\Pi^{-1}(q)$ for each $q\in  X'$. Let $a^i_j\in k$ be such that
\[
\nabla_ve_j=-a^i_je_i
\]
(recall Remark \ref{einsconv}). Thus, in this frame, the equation $\nabla_vX=0$ is equivalent to $X'=AX$, where $A=(a^i_j)$ and $X'=v(x^i)e_i$ if $X=x^ie_i$, $x^i\in k$.

We review the construction of a Picard-Vessiot extension. We define the differential ring extension $(k[X^i_j,\frac{1}{\det}],\tilde{v})$ of $(k,v)$, where $(X^i_j)$ is an $n\times n$ matrix of indeterminates, $\det:=\det(X^i_j)$ and
\[
\tilde{v}(X^i_j)=a^i_kX^k_j.
\]
Note that we can make $GL_n(\mathbb{C})$ act on $k[X^i_j,\frac{1}{\det}]$ through differential automorphisms over $(k,v)$ by setting for $(g^i_j)\in GL_n(\mathbb{C})$
\begin{eqnarray*}
(g^i_j): k[X^i_j,\frac{1}{\det}] & \longrightarrow & k[X^i_j,\frac{1}{\det}]\\
X^i_j   & \longmapsto & X^i_lg^l_j
\end{eqnarray*}
A Picard-Vessiot extension of $k$ for the matrix differential equation $X'=AX$ is given by the quotient field of
\[
k[X^i_j,\frac{1}{\det}]/ I,
\]
where $I$ is a maximal differential ideal. Since $GL_n(\mathbb{C})$ acts through differential automorphisms, the action permutes the maximal differential ideals. A representation of the differential Galois group of $\nabla$ is given by~\cite[Corollary 4.10]{magid}
\[
G=\{(g^i_j)\in GL_n(\mathbb{C})|\ (g^i_j):I\mapsto I\}
\]
the stabilizer of $I$ under this action.

Let us assume that the holomorphic chart $(U,z)$ and $v$ are such that $v$ restricted to $U$ coincides with $\frac{\partial}{\partial z}$. Fix an injection $\iota: k\rightarrow\mathbb{C}[\frac{1}{z}][[z]]$. Identifying $v^i_j(z)$ with their power series expansion we have a differential homomorphism
\begin{eqnarray*}
\Phi:(k[X^i_j,\frac{1}{\det}],\tilde{v}) & \longrightarrow & (\mathbb{C}[\frac{1}{z}][[z]],\frac{\partial}{\partial z})\\
X^i_j   &\longmapsto & v^i_j(z)
\end{eqnarray*}
such that $\Phi(f)=\iota(f)$ if $f\in k$, i.e.
\[
\xymatrix{
k[X^i_j,\frac{1}{\det}]\ar[rr]^\Phi & & \mathbb{C}[\frac{1}{z}][[z]]\\
k\ar[u]\ar[urr]_\iota & &
}
\]
If we set $I=\ker(\Phi)$, then $I$ is a maximal differential ideal and the choice of $v^i_j(z)$ induces the representation of the Galois group of $\nabla$ by $G$.

In order to state our theorem we need to introduce the following concepts~\cite{We}~\cite{ber}:

\begin{defn}
Let $P[X^i_j]\in\mathbb{C}[X^i_j]$ be a homogeneous non-constant polynomial. If $P[v^i_j](z)=\iota(f)$ for some $f\in k$ we say that $f$ is a \emph{dual first integral} of $\nabla$ with degree defined to be the degree of $P$. We denote by $k_\nabla$ the field generated over $\mathbb{C}$ by the quotients of dual first integrals of the same degree.
\end{defn}

\begin{defn}\label{def0}
The connection $\nabla$ is called:
\begin{itemize}
\item \emph{standard} if $k(v^i_j)=\mathbb{C}(v^i_j)$ and $k$ is Galois over $k_\nabla$.
\item \emph{basic} if $k=k_\nabla$.
\end{itemize}
\end{defn}

\begin{rem}
Since $k_\nabla\subseteq\mathbb{C}(v^i_j)$ then if $\nabla$ is basic we have $\mathbb{C}(v^i_j)=k_\nabla(v^i_j)=k(v^i_j)$. In particular basic connections are standard.
\end{rem}

\begin{rem}
As a corollary of Lemma \ref{lem1b}, we will see that, under the unimodular and reductive hypothesis, the automorphisms of $k$ over $k_\nabla$ can be identified as a subgroup of the symmetries of $\nabla$ provided $k$ is standard. In particular, this means that $\nabla$ would descend to a connection over $k_\nabla$. Note that $k$ and $k_\nabla$ are $\mathbb{C}$-algebras. Another approach to the descent problem is seen in ~\cite{vHoeijvdP}, where the treatment is in terms of Galois cohomology and deals with descent on the field of constants.
\end{rem}

\begin{thm}\label{theorem}
Let $\nabla$ be a standard connection. If the Galois group of $\nabla$ is unimodular and reductive the sequence
{\small
\[
1\longrightarrow Z(G)\longrightarrow G\longrightarrow Aut_{G_F}(G)\longrightarrow
\textrm{Aut}_\nabla( X)\longrightarrow 1
\]}
is exact. 
\end{thm}

\begin{rem}
In the statement above $Aut_{G_F}(G)$ denotes the image, in the group of automorphisms of $G$, of the normalizer $N_{G_F}(G)$ of $G$ in $G_F$.  $Z(G)$ denotes the center of $G$.
\end{rem}

The remainer of the monograph is devoted to a proof of this theorem. The hypothesis on the differential Galois group allows us to use the following result ~\cite{compoint}:

\begin{thm}[Compoint]\label{compthmb}
When $G$ is reductive and unimodular the ideal $I$ is generated by the $G$-invariants contained therein. Moreover, if $P_1,\ldots, P_r$ is a set of homogeneous generators for the $\mathbb{C}$-algebra of $G$-invariants, with respective degrees $n_1,\ldots,n_r$, in $\mathbb{C}[X^i_j]$, and if $f_1,\ldots, f_r\in k$ are such that $P_i-f_i\in I$, then $I$ is generated over  $k[X^i_j,\frac{1}{\det}]$ by $P_i-f_i$, where $i\in\{1,\ldots, r\}$.
\end{thm}

\begin{rem}
In~\cite{compoint} and in~\cite{beukers} the statement of this theorem is reserved for $k=\mathbb{C}(z)$ and $v=\frac{\partial}{\partial z}$. But the proof by F. Beukers in ~\cite{beukers} carries through \emph{mutatis mutandis} when $\mathbb{C}(z)$ is replaced by $k$. The careful reader would note that F. Beukers does not make the unimodular hypothesis explicit in his paper, but the assumption is needed to guarantee the invariance of $\det$.
\end{rem}

\begin{rem}\label{rem0}
In the notation of the theorem we obtain
\[
k_\nabla=\mathbb{C}({\footnotesize \frac{f_i^{m_i}}{f_j^{m_j}}})_{\{(i,j)\in\{1,\ldots,r\}^2|\ m_in_i=m_jn_j,\ f_j\ne 0\}}.
\]
\end{rem}

\begin{prop}\label{lem0b}
If $\nabla$ is such that $k(v^i_j)=\mathbb{C}(v^i_j)$ and $1$ is a dual first integral of $\nabla$ (i.e. if there is a homogeneous $P(X^i_j)\in\mathbb{C}[X^i_j]^G\setminus\mathbb{C}$ such that $P(v^i_j)=1$), then $k/k_\nabla$ is Galois with abelian Galois group. In particular, $\nabla$ is standard.
\end{prop}

\noindent\emph{Proof}:
The assumption on $\nabla$ implies that $\mathbb{C}(v^i_j)^G=k(v^i_j)^G=k$. In the notation of Theorem \ref{compthmb} we therefore have $\mathbb{C}(v^i_j)^G=\mathbb{C}(f_1,\ldots,f_r)$. By hypothesis there is a homogeneous $P(X^i_j)\in\mathbb{C}[X^i_j]^G$ such that $P(v^i_j)=1$. Denote by $d$ the degree of $P(X^i_j)$, and by $[d,n_l]:=dm_l=n_ld_l$ the least common multiple of $d$ and of the degree of $P_l(X^i_j)$, $n_l$. Now
\[
f_l^{d_l}=\frac{f_l^{d_l}}{1^{m_l}}\in k_\nabla
\]
Thus $k$ is the splitting field over the $\mathbb{C}$-algebra $k_\nabla$ of the polynomial 
\[
\prod_{l=1}^r (Z^{d_l}-f_l^{d_l}).
\]
\hfill$\bigstar$

\begin{rem}
It is easy to produce examples to justify the unimodular hypothesis of the theorem. Indeed, consider
\[
\left(\begin{array}{c} y_1 \\ y_2\end{array}\right)' =
\left(\begin{array}{cc}2z & 0 \\
                       0 & \frac{1}{z}+2z\end{array}\right)
\left(\begin{array}{c} y_1 \\ y_2\end{array}\right).
\]
A fundamental system of solutions is given by
\[
\left(\begin{array}{cc}e^{z^2} & 0 \\
                       0 & ze^{z^2}\end{array}\right),
\]
and the connection is therefore basic. Indeed $Y^2_2/Y^1_1=z$. The Galois group is $\mathbb{G}_{m,\mathbb{C}}$, which is represented by
\[
\left\{\left(\begin{array}{cc}\lambda & 0\\ 0 & \lambda\end{array}\right)\Big|\ \lambda\in\mathbb{C}^*\right\}.
\]
The ideal of homogenous polynomials vanishing at the solutions is generated by $zX^1_1-X^2_2$ , $X^1_2$ and $X^2_1$. There is no vanishing $G$-invariant homogeneous polynomial of order greater than $0$. It follows that $G_F$ is $GL_2(\mathbb{C})$, $N_{G_F}(G)=G_F$, $G_F/G=PGL_2(\mathbb{C})$, but the connection is symmetric only with respect to $z\mapsto -z$. Indeed, one can lift this symmetry by mapping $(y_1,y_2)$ into $(y_1,-y_2)$.
\end{rem}

\begin{rem}\label{exa0}
An example justifying the requirement on $\nabla$ to be standard is a little more delicate; when one lifts a standard equation with symmetries by means of a Galois covering, the symmetries are lifted together with the equation. So the way to obtain the example is by lifting an equation with symmetries, say the standard $D_{2\cdot 4}$ equation ~\cite{ber} defined on the $x$-sphere, through a non-Galois covering, say through $x\mapsto z(z-1)^2$, to the $z$-sphere.
\end{rem}

\section{Geometric construction: a covering space with covering group $Aut_{G_F}(G)$}

\begin{rem}
We keep the same notation from the previous section, we will assume, from now on, that $\nabla$ has reductive
unimodular Galois group and that $\nabla$ is standard or basic.
\end{rem}

Following J.A. Weil ~\cite{weil} we will work with first integrals (i.e. with the solutions of the adjoint
system) rather than with horizontal sections. Recall that a linear first integral
of $\nabla'_{\frac{\partial}{\partial z}}$ over $\Pi^{-1}U$ refers to a
function $\Phi:\Pi^{-1}U\to \mathbb{C}$ of the form
\[
\Phi(z,y^1,\ldots,y^n)=\frac{\partial \Phi}{\partial y^j}(z)y^j,
\]
where $\frac{\partial \Phi}{\partial y^j}(z)$ is holomorphic in $U$ and constant on the horizontal sections of  $\nabla'_{\frac{\partial}{\partial
z}}$, i.e. if $\nabla'_{\frac{\partial}{\partial z}}X=0$ then $z\mapsto \Phi(X(z))$ is a constant function.

Let $\Phi^i$, $i\in\{1,\dots,n\}$, be a system of linear first integrals of
$\nabla'_{\frac{\partial}{\partial z}}$ over $\Pi^{-1}U$ such that
$\frac{\partial\Phi^i}{\partial y^j}$ is an invertible matrix. We call such a system a \emph{fundamental system} of linear first integrals.

Let $F$ be the homogeneous ideal in $\mathbb{C}[X^i_j,\frac{1}{\det}]^G$ generated by the homogeneous polynomials of $\mathbb{C}[X^i_j]^G$
vanishing at $\frac{\partial \Phi^i}{\partial y^j}$.
\begin{lem}
The ideal $F$ is prime.
\end{lem}

\noindent\emph{Proof}:
Without loss of generality we may assume that $\frac{\partial \Phi^1}{\partial y^1}\ne 0$, so that $X^1_1$ is not in the kernel of the evaluation map $\mathbb{C}[X^i_j,\frac{1}{\det}]\rightarrow k(v^i_j): X^i_j\mapsto \frac{\partial \Phi^i}{\partial y^j}$. So we define the map
\begin{eqnarray*}
\mathbb{C}[\frac{X^i_j}{X^1_1},\frac{(X^1_1)^n}{\det}] & \longrightarrow & k(\frac{\partial \Phi^i}{\partial y^j})\\
                                   \frac{X^i_j}{X^1_1} & \longmapsto     & \frac{\frac{\partial \Phi^i}{\partial y^j}}{\frac{\partial \Phi^1}{\partial y^1}}.
\end{eqnarray*}
Every element in the kernel of this map defines, after homogenizing with $X^1_1$, a homogeneous element in the kernel of the evaluation map. Conversely, every homogeneous element of degree $d$ in the kernel of the evaluation map, after dividing by $(X^1_1)^d$, defines an element in the kernel of this map. Thus, since $k(\frac{\partial \Phi^i}{\partial y^j})$ is an integral domain and the homogeneous elements in the kernel of the evaluation map form a prime ideal. The lemma now follows by letting $G$ act and taking invariants.\hfill$\bigstar$

Note that if $AX=X'$ is the matrix equation form of $\nabla_vX=0$, the adjoint system\footnote{The common notation is $X'=-A^tX$ but our notation is aimed to make more transparent the interplay between both systems.} for the linear first integrals is given by $X(-A)=X'$. Indeed, if $S$ is a fundamental system of solutions for $AX=X'$ then from
\[
S^{-1}S=I_n
\]
and we see that $S^{-1}$ is a fundamental system of linear first integrals. Moreover
\[
0=I_n'=(S^{-1}S)'=(S^{-1})'S+S^{-1}S'=(S^{-1})'S+S^{-1}AS
\]
and therefore $(S^{-1})'=-S^{-1}A$. When considering first integrals we use the differential ring extension $(k[X^i_j,\frac{1}{\det}],\overline{v})$ of $(k,v)$, where $\overline{v}X^i_j=-X^i_ka^k_j$, and we let $GL_n(\mathbb{C})$ act to the left.

A maximal differential ideal $I$ determines a representation of the Galois group $G\subseteq GL_n(\mathbb{C})$ given by the stabilizer of $I$. On the other hand $I$ also determines the homogeneous ideal $F\subseteq\mathbb{C}[X^i_j]^G$ generated by the homogeneous elements in $I\cap\mathbb{C}[X^i_j]^G$. According to Compoint's Theorem, any maximal differential ideal of the differential ring $(k[X^i_j,\frac{1}{\det}],\overline{v})$, with stabilizer $G$ and defining the ideal $F$, is uniquely determined by a $\mathbb{C}$-homomorphism $\phi: R\rightarrow k$. Indeed, in the notation of the statement of Compoint's theorem we have a map
\begin{eqnarray*}
\mathbb{C}[X^i_j,\frac{1}{\det}]^G & \longrightarrow & k\\
           P_i                     & \longmapsto     & f_i.
\end{eqnarray*}
The ideal $F$, by definition, is contained in the kernel, so this map factors through a unique map $\phi: R\rightarrow k$.

\begin{lem}\label{lemact}
The group $H:=Aut_{G_F}(G)$ acts on $k_\nabla$.
\end{lem}

\noindent\emph{Proof}:
The group $H=N_{G_F}(G)/Z(N_{G_F}(G))$ acts on the Zariski open subset of $\textrm{Proj}(\mathbb{C}[X^i_j])$ defined by $\det\ne 0$, and by the definition of $G_F$ this action fixes the Zariski closed subset defined by the homogeneous ideal generated by $F$ in $\mathbb{C}[X^i_j]$. So passing to the quotient by the action of $G/Z(G)$, we have that $H$ acts on $\textrm{Proj}(\mathbb{C}[X^i_j]^G)$, leaving the variety defined by the homogeneous prime ideal $F$ invariant. In particular, it defines an automorphism of the field of meromorphic functions over this variety, which is actually $k_\nabla$.\hfill$\bigstar$

In particular, if $\nabla$ is basic the group also describes automorphisms of $k$ (cf. Definition \ref{def0}).

\[
\xymatrix{
k[X^i_j,\frac{1}{\det}]/I \\
k\ar@{-}[u] \\
k_\nabla\ar@{-}[u] \\
(k_\nabla)^H\ar@{-}[u]
}
\]

We will now give a geometric construction that will make clear the interplay between all the groups and fields involved.

Fix a meromorphic tangent vector field $v$ over $ X$. Let $p\in X$ over which $\nabla_v$ is not singular and fix a frame of holomorphic horizontal sections $V_1,\ldots,V_n$ of $\nabla_v$ around $p$. We lift the connection to the frame bundle $GL(\Pi)$.

We extend $U\subseteq X$ to the maximal open set $ X'$ over which the frame $(V_1,\ldots,V_n)$ can be extended holomorphically (as a multi-valued frame).

The extension of $(V_1,\ldots,V_n)$, together with its orbit under the Galois group $G$, defines a subsheaf $\mathfrak{F}$ of $GL(\Pi)\upharpoonright_{ X'}$ whose sections are horizontal holomorphic frames under $\nabla_v$. The sheaf $\mathfrak{F}$ gives rise to a regular covering of $ X'$ with covering group $G$ ~\cite{malgrange}~\cite[Th\'eor\`eme 5.3.1]{casale}. Denote by $\widetilde{ X'}$ the covering space corresponding to the center $Z(G)$ of $G$.\\
The diagram above implies the following tower of covering spaces:
\[
\xymatrix{
\widetilde{ X'} \\
 X'\ar@{-}[u] \\
 X'_\nabla\ar@{-}[u] \\
( X'_\nabla)^H\ar@{-}[u]
}
\]
where $ X'_\nabla$ corresponds to the projection of $ X'$ to the Riemann surface with meromorphic functions $k_\nabla$. In the case that $\nabla$ is basic we have  $ X'= X'_\nabla$. So we obtain $\widetilde{ X'}$ covering $ X'$ with covering group $PG$, and $\widetilde{ X'}$ covering $( X'_\nabla)^H$ with covering transformations $H=Aut_{G_F}(G)$. By the Galois correspondence, if $\nabla$ is basic, the covering group of $ X'$ over $( X'_\nabla)^H$ is given by the outer-automorphisms of $G$ in $G_F$. If $\nabla$ is standard then the outer-automorphisms of $G$ in $G_F$ define covering transformations of $ X'$ over $( X'_\nabla)^H$ since $G$ fixes $ X'_\nabla$ and $ X'$ is Galois over $ X'_\nabla$.

\section{The map $Aut_{G_F}(G)\longrightarrow Aut_\nabla( X)$}\label{map}

Given a maximal differential ideal $I\subseteq k[X^i_j,\frac{1}{\det}]$, we obtain a representation of the Galois group, namely,  the stabilizer $G\subseteq GL_n(\mathbb{C})$ of $I$. Any $g\in GL_n(\mathbb{C})$ sends $I$ into another maximal differential ideal $I^g:=g(I)$. The stabilizer of $I^g$ is now the conjugate $gGg^{-1}$ of $G$. So the collection of maximal differential ideals whose stabilizer is $G$, or equivalently, whose representation of the Galois group is $G$, will be given by the orbit of $I$ under $N_{GL_n(\mathbb{C})}(G)$. We set 
\[
\mathcal{I}:=\{I^g|\ g\in N_{GL_n(\mathbb{C})}(G)\}.
\]
On the other hand, each maximal differential ideal $I^g\in\mathcal{I}$ defines the homogeneous ideal $F^g:=g(F)\subseteq \mathbb{C}[X^i_j]^G$ given by the homogeneous elements in $I^g\cap \mathbb{C}[X^i_j]^G$. So we consider the sub-collection
\[
\mathcal{F}:=\{I^g\in\mathcal{I}|\ F=F^g\}
\] 
of $\mathcal{I}$ which is precisely the orbit of $I$ under the elements of $G_F\cap N_{GL_n(\mathbb{C})}(G)$. So we have
\[
\mathcal{F}:=\{I^g |\ g\in N_{G_F}(G)\}.
\]
Compoint's Theorem implies that for every $I^g\in\mathcal{F}$ we have a unique map
\[
\phi_g: R\rightarrow k.
\]
Note that since $g\in N_{G_F}(G)$ fixes $F\subseteq R$, it defines an automorphism of $R$. In particular,
\[
\phi\circ g=\phi_g.
\]

\begin{rem}
The open set $ X'$ from our previous section is the variety given by $ X\setminus S$, where $S$ is the collection of points where $\phi(\det)$ vanishes, together with the singular points of $\nabla$.
\end{rem}

Consider a $g=(g^i_j)\in N_{G_F}(G)$,  so that now we have two $k$-valued points $\phi$ and $\phi_{g^{-1}}$. The $k$-valued point $\phi=\phi_e$ corresponds to the ideal $I$ and $\phi_{g^{-1}}$ to $I^{g^{-1}}$. Under the action notation we write $\phi=\phi_{g^{-1}}\circ(g^i_j)$. Similarly, one can see $(g^i_j)$ acting on $k_\nabla$ (Lemma \ref{lemact}). Indeed, instead of seeing $(g^i_j)$ as acting linearly on $G$-invariant polynomials, one can see it as acting linearly on quotients of rational first integrals of same degree. On the other hand, $k$ is Galois over $k_\nabla$, so the automorphism $(g^i_j)$ on $k_\nabla$ lifts to an automorphism of $k$. We denote by $\sigma^*$ such a lifting $(g^i_j)$. In this fashion we obtain a commutative diagram:
\[
\xymatrix{
R\ar[rr]^{(g^i_j)}\ar[d]_{\phi_{g^{-1}}}\ar[drr]^{\phi} & & R\ar[d]^{\phi_{g^{-1}}}\\
k\ar[rr]_{\sigma^*} & & k
}
\]
Reversing the arrows:
\[
\xymatrix{
 X'_\nabla\ar@{<-}[rr]^{(g^i_j)^*}\ar@{<-}[d]_{\phi_{g^{-1}}^*}\ar@{<-}[drr]^{\phi^*} & &  X'_\nabla\ar@{<-}[d]^{\phi_{g^{-1}}^*}\\
 X'\ar@{<-}[rr]_{\sigma} & &  X'
}
\]

\begin{lem}\label{lem1b}
The automorphism $\sigma$ (or equivalently a lifting of $(g^i_j)^*$) is a symmetry of $\nabla$.
\end{lem}

\noindent\emph{Proof}:
If $P[X^i_j]\in\mathbb{C}[X^i_j]^G$ and $f\in k$ are such that $\phi_{g^{-1}}(P[X^i_j])=f$, then
\begin{eqnarray*}
\phi(P[X^i_j]) & = & \phi_{g^{-1}}(P[g^i_kX^k_j])  \\
               & = & \sigma^*\phi_{g^{-1}}(P[X^i_j])
\end{eqnarray*}
i.e.,
\begin{eqnarray*}
\phi(P[X^i_j]) & = & \sigma^*f.
\end{eqnarray*}

We again take holomorphic charts $(U,z)$ and $(\pi^{-1}U,z,y^1,\ldots,y^n)$ as before in Section \ref{secdef}. Set $V=\sigma(U)$ and consider a chart of $E$ giving a vector bundle trivialization $(\Pi^{-1}V,w,x^1,\ldots,x^n)$ with $w\circ \sigma = z$. Note that
\[
\sigma_*\frac{\partial}{\partial z}(w)=\frac{\partial}{\partial z}(w\circ\sigma)=\frac{\partial}{\partial z}(z)=1,
\]
hence $\sigma_*\frac{\partial}{\partial z}=\frac{\partial}{\partial w}$.

Let $p\in U$ ($\sigma(p)\in V$) and consider some linear first integrals $\frac{\partial \Psi^i}{\partial x^j}$ defined by $\phi_{g^{-1}}$ in $\Pi^{-1}V$. Since the $k$-valued point given by Compoint's Theorem is nothing other than the restriction of the evaluation homomorphism, the equalities above implies
\begin{eqnarray*}
\frac{\partial \Phi^i}{\partial y^j} (p) & = & g^i_k\frac{\partial \Psi^k}{\partial x^j} (p)\\
                                         & = & \frac{\partial \Psi^i}{\partial x^j}(\sigma(p)),
\end{eqnarray*}
hence that
\begin{eqnarray}\label{fibmor}
\frac{\partial \Phi^i}{\partial y^j}(p) & = & \frac{\partial \Psi^i}{\partial x^j}(\sigma(p)) = \frac{\partial \Psi^i}{\partial y^l}(\sigma(p))\frac{\partial y^l}{\partial x^j}(\sigma (p)).
\end{eqnarray}
But $\frac{\partial \Phi^i}{\partial y^j}$ is invertible in $U$, whence
\begin{eqnarray*}
\frac{\partial x^j}{\partial y^l}(\sigma(p)) & = & \frac{\partial y^j}{\partial \Phi^i}(p)\frac{\partial \Psi^i}{\partial y^l}(\sigma(p))\\
                                             & = & f^j_l(p).
\end{eqnarray*}
Letting $G$ act on both sides of (\ref{fibmor}), it follows that $f^j_l$ describes (by analytic extension) a meromorphic function over $ X$ (Galois correspondence). So the map $U\rightarrow V$ described by
\begin{eqnarray*}
w(z,y^1,\ldots,y^n) & = &\sigma(z)\\
x^j(z,y^1,\ldots,y^n) & = & f^j_i(z)y^i
\end{eqnarray*}
gives the transform on the fiber coordinates that lifts the covering transformation $\sigma$ to a horizontal automorphism.\hfill$\bigstar$

If there is another $\sigma'$ with the property $\phi^*=\phi^*_{g^{-1}}\circ\sigma'$ then
\[
\phi^*\circ\sigma^{-1}\circ\sigma'=\phi^*_{g^{-1}}\circ\sigma\circ\sigma^{-1}\circ\sigma'=
\phi^*_{g^{-1}}\circ\sigma'=\phi^*.
\]
This says that $\sigma^{-1}\circ\sigma'$ is a covering transformation of $\phi^*$, and so $\sigma$ and $\sigma'$ are two liftings of the automorphism of $k_\nabla$ defined by $(g^i_j)$. The previous lemma asserts that $\sigma'$ is a symmetry. Conversely, given a symmetry of the form in the previous lemma (e.g. the identity or $\sigma$), then for any covering transformation of $\phi^*$ (or equivalently any Galois automorphism of $k$ over $k_\nabla$), the symmetry composed with the covering transformation gives us another symmetry. So we conclude that $\nabla$ descends to a connection over the compact Riemann surface with field of meromorphic functions $k_\nabla$. This last connection with Picard-Vessiot extension $k_\nabla(v^i_j)=\mathbb{C}(v^i_j)$ is by construction a basic connection. We obtain the following proposition:

\begin{prop}
Any standard connection $\nabla$ is the pullback of a basic connection over the Riemann surface with field of meromorphic functions $k_\nabla$.
\end{prop}

Finally, as $H=Aut_{G_F}(G)$ acts on $k_\nabla$ as symmetries of this last basic connection (Lemma \ref{lem1b}), we can in turn descend the vector bundle and connection all the way down to the Riemann surface $ X_0$ with field of meromorphic functions $(k_\nabla)^H=:k_0$. Denote by $E_0$ this quotient vector bundle over $ X_0$ and by $\nabla_0$ the resulting quotient connection. Then by definition $(E,\nabla)$ is the pullback of $(E_0,\nabla_0)$, and we have the following tower of Galois extensions:
\[
\xymatrix{
k(v^i_j)^{Z(G)}\\ k\ar@{-}[u]\\ k_0\ar@{-}[u]
}
\]

\begin{lem}
$k_0(v^i_j)=k(v^i_j)$ is a Picard-Vessiot extension for $\nabla_0$. The Galois group of $\nabla_0$ is represented by a subgroup of $N_{G_F}(G)$ with projective Galois group $H=Aut_{G_F}(G)$.
\end{lem}

\noindent\emph{Proof}:
The first statement on the lemma follows from the fact that $\nabla$ is the pullback of $\nabla_0$. The establish the second assertion it suffices to notice $k_0=(k_\nabla)^H=(\mathbb{C}(v^i_j)^{Z(G)})^H$.\hfill$\bigstar$

To obtain the map $Aut_{G_F}(G)\rightarrow Aut_\nabla( X)$ identify the group of symmetries with the Galois group $Aut_{k_0}(k)$, where in the map is seen to arise from the Galois action on the tower immediately before the statement of the lemma.

\section{Right exactness}

Take $\sigma\in Aut_\nabla( X)$. Put $V=\sigma(U)$ and consider a chart of $E$ giving a vector bundle trivialization $(\Pi^{-1}V,w,x^1,\ldots,x^n)$ with $w\circ \sigma = z$. Denote by $\nabla'$ the holomorphic connection on $V$ induced by $\nabla$. Since $\sigma$ is a symmetry there is a horizontal vector bundle isomorphism
\begin{eqnarray*}
\tilde{\sigma}:(\Pi^{-1}U,\nabla') & \longrightarrow & (\Pi^{-1}V,\nabla')
\end{eqnarray*}
lifting $\sigma$. Note that $\tilde{\sigma}^*w=w\circ \tilde{\sigma}= (w\circ\Pi) \tilde{\sigma}= w\circ\sigma\Pi = z\circ\Pi=z$ and $\sigma_*\frac{\partial}{\partial z}=\frac{\partial}{\partial w}$.

The hypothesis $\sigma\in Aut_\nabla( X)$ implies that each pullback $\tilde{\sigma}^*\Psi^i$ along $\tilde{\sigma}$ of a fundamental system of linear first integrals
$\Psi^i$ of $\nabla'_{\frac{\partial}{\partial w}}$ over $\Pi^{-1}V$ is a linear first integral of $\nabla'_{\frac{\partial}{\partial z}}$ over $\Pi^{-1}U$. In other words, there exist $c^i_k$ such that
\begin{eqnarray*}
c^i_k\Phi^k(z,y^1,\ldots,y^n) & = & \tilde{\sigma}^*\Psi^i(z,y^1,\ldots,y^n)\\
c^i_k\frac{\partial \Phi^k}{\partial y^j}(z)y^j& = & \frac{\partial \Psi^i}{\partial x^j}(\sigma(z))\tilde{\sigma}^*x^j(y^1,\ldots,y^n).
\end{eqnarray*}
Taking $x^i$ such that $\tilde{\sigma}^*x^j(y^1,\ldots,y^n)=y^j$ we have
\begin{eqnarray}\label{pre}
c^i_k\frac{\partial \Phi^k}{\partial y^j}(z)& = & \frac{\partial \Psi^i}{\partial x^j}(\sigma(z)).
\end{eqnarray}

Choose $f\in k$ such that $f(\sigma (q))v_{\sigma (q)}=(\sigma_*v)_{\sigma (q)}$, for every $q\in X$. If $AX=X'$ is the matrix equation form of $\nabla_vX=0$, then on $U$ we have
\[
\nabla' X=(AX-X')\otimes dz
\]
Let $\gamma$ be a path from $p\in U$ to $\sigma (p)\in V$ avoiding the singularities of $\nabla$. If $\nabla'_{\frac{\partial}{\partial z}}X=0$, and if $X$ is analytically extended along $\gamma$, we have
\begin{eqnarray*}
\nabla'_{\sigma_*\frac{\partial}{\partial z}}X & = & \langle (AX-X')\otimes dz, \sigma_*\frac{\partial}{\partial z} \rangle\\
  & = & \langle (AX-X')\otimes dz, f(z)\frac{\partial}{\partial z} \rangle\\
  & = & f(z)\nabla'_{\frac{\partial}{\partial z}}X\\
  & = & 0.
\end{eqnarray*}
So we may take $\Psi^i$ as the analytic extension of $\Phi^i$ along $\gamma$, and (\ref{pre}) becomes:
\begin{eqnarray}\label{post}
c^i_k\frac{\partial \Phi^k}{\partial y^j}(z)& = & \frac{\partial \Phi^i}{\partial x^j}(\sigma(z))
\end{eqnarray}

\begin{lem}
The matrix $(c^i_k)$ defined on (\ref{post}) is in the normalizer of $G$ in $G_F$.
\end{lem}

\noindent\emph{Proof}:
As in the previous section, let $ X_0$ be the Riemann surface obtained as the quotient space of $ X$ by the group of symmetries $Aut_\nabla( X)$. By the definition of symmetry there exists a vector bundle $E_0$ and a connection $\nabla_0$ such that $(E,\nabla)$ is the pullback of $(E_0,\nabla_0)$ under the covering map induced by the action.

The projection $\gamma_0$ of $\gamma$ to $ X_0$ is a closed curve, and $(c^i_k)$ therefore defines an element of the monodromy of  $\nabla_0$.

Let $P[X^i_j]\in\mathbb{C}[X^i_j]^G$, i.e. a $G$-invariant polynomial such that
\[
P[\frac{\partial \Phi^i}{\partial y^j}](z)=\iota(f)(z)
\]
for some $f\in k$.

If $f=0$, then under analytic continuation along $\gamma_0$ it remains the case that $P[c^i_k\frac{\partial \Phi^k}{\partial y^j}(z)]=0$. This implies $(c^i_k)\in G_F$.

On the other hand, $k$ is an extension of the field of meromorphic functions over $ X_0$. Under analytic continuation along $\gamma_0$, an arbitrary non-zero $f$ is mapped into $f_\sigma$ by a covering (i.e. Galois) automorphism of $ X$ (of $k$) over $ X_0$ (over $\mathbb{C}( X_0)$). Whence $P[c^i_k\frac{\partial \Phi^k}{\partial y^j}(z)]=\iota(f_\sigma)(z)$, and so by the Galois correspondence $P[X^i_j]$ is invariant under $G^{(c^i_k)}$ (the conjugate of $G$ by $(c^i_k)$).

A symmetric argument allows us to conclude that the invariant polynomials under $G$ and under $G^{(c^i_k)}$ coincide. So Compoint's Theorem implies $G^{(c^i_k)}=G$. This completes the proof.\hfill$\bigstar$

The Lemma implies the map defined in Section \ref{map} is surjective, and the theorem follows.

\section{Examples}

\begin{exa}\label{ex1}
Consider the elliptic curve defined by the equation $z^3-z=w^2$ with field of meromorphic functions $\mathbb{C}(z,w)$. We take its invariant differential form $\frac{dz}{2w}$ and its dual tangent vector field $v=2w\frac{\partial}{\partial z}$. We denote by $\sigma$ the fourth order automorphism $\sigma: w\mapsto iw, z\mapsto -z$; so that $\sigma_*v=-iv$.

The differential equation

{\scriptsize
\[
v^3(y)-\frac{3(5z^2-3)}{w}v^2(y)+\frac{1}{3}\frac{194z^4-230z^2+108}{w^2}v(y)-\frac{4}{27}\frac{364z^6-665z^4+1030z^2-405}{w^3}y=0
\]}
has differential Galois group $G_{27}$. Indeed, written in terms of $z$ and $\frac{\partial}{\partial z}$, the equation corresponds to:
{\small
\[
y'''-\frac{3}{z}y''+\frac{1}{12}\frac{77z^4-122z^2+81}{(z^3-z)^2}y'-\frac{1}{54}\frac{364z^6-665z^4+1030z^2-405}{(z^3-z)^3}y=0.
\]
}

This equation is irreducible and has unimodular Galois group. Furthermore to get the Galois group using methods of \cite{ulmer} we notice that it has a two-dimensional space of third degree invariants, which corresponds to the dual first integrals $0$ and $z^3w^3$. The wronskian of the equation in terms of $v,z,w$ is $zw$, the ratio of two sixth degree invariants is $z$ and the ratio of two ninth degree invariants is $w$, so the induced connection is basic.

If $X$, $Y$ and $Z$ are the solution based at point $p:(z(p),w(p))=(0,0)$, with respectively $z$-exponents $2$, $3/2$ and $5/2$, and leading coefficient $1$, then
\[
Y^2Z+X^2Y-\frac{1}{81}Z^3=0.
\]
This vanishing homogenous polynomial is a $G_{27}$ invariant which describes an elliptic curve, and the subgroup of $GL_3({\mathbb{C}})$ leaving this elliptic curve invariant is the lifting of $F_{36}\subset PSL_3(\mathbb{C})$. So the group $G_F/Z(G_F)$ is $F_{36}$. The group $PG_{27}$ (the projective version of $G_{27}$) is a normal subgroup of index four in $F_{36}$; the quotient is a cyclic group of order four. On the other hand the equation is invariant under $\sigma$ and we obtain
{\small
\[
1\longrightarrow Z(G_{27}) \longrightarrow G_{27}\longrightarrow F_{36}\longrightarrow \langle \sigma \rangle\longrightarrow 1.
\]
}
The equation descends into the Riemann sphere parameterized by $x=z^2$ to an equation with Galois group $F_{36}^{SL_3}$, i.e.
\[
y'''+\frac{1}{48}\frac{41x^2-50x+45}{x^2(x-1)^2}y'-\frac{1}{432}\frac{364x^3-665x^2+1030x-405}{x^3(x-1)^3}y=0,
\]
which is given in terms of $x$ and $\frac{\partial}{\partial x}$.

Algebraically we have $v=2w\frac{\partial}{\partial z}=4wz\frac{\partial}{\partial x}$ and $\mathbb{C}(z,w)\supseteq\mathbb{C}(z)\supseteq\mathbb{C}(x)$, with $\mathbb{C}(z)=\mathbb{C}(z,w)^{\sigma^2}$ and $\mathbb{C}(x)=\mathbb{C}(z,w)^{\sigma}$. We also have that $\frac{\partial}{\partial z}$ is invariant under $\sigma^2$ and that $\frac{\partial}{\partial x}$ is invariant under $\sigma$. We thus have the following tower of linear operators:
\[
\xymatrix{
\mathbb{C}(z,w)[v]\ar@{-}[d] & \\
\mathbb{C}(z)[\frac{\partial}{\partial z}]\ar@{-}[d] & =\left(\mathbb{C}(z,w)[v]\right)^{\sigma^2}\\
\mathbb{C}(x)[\frac{\partial}{\partial x}] & =\left(\mathbb{C}(z,w)[v]\right)^{\sigma}
}
\]
The first of the equations corresponds to an operator in $\mathbb{C}(w,z)[v]$, which when written in terms of $\frac{\partial}{\partial z}$ defines the operator in the second equation, which is in $\mathbb{C}(z)[\frac{\partial}{\partial z}]$. Finally writing it in terms of $\frac{\partial}{\partial x}$, we obtain the operator in $\mathbb{C}(x)[\frac{\partial}{\partial x}]$ corresponding to the third equation. So the operator in the first equation is actually an operator in the bottom of the tower.
\end{exa}

\begin{exa}
Consider the differential equation $L(y)=0$ given by:
\[
y''-\frac{z^4-3z^2-1}{1+z^4}y=0.
\]
Two linearly independent solutions are given by:
\[
Y_1=\sqrt[4]{z^4-1}e^{\int\frac{1}{\sqrt{z^4-1}}},\quad
Y_2=\sqrt[4]{z^4-1}e^{-\int\frac{1}{\sqrt{z^4-1}}}.
\]
So according to the Kovacic algorithm, the Differential Galois group is the infinite dihedral group $D_\infty$. A basis of invariants of this group on $\mathbb{C}[X^i_j]$ is given by
\begin{eqnarray*}
X_{2,1}:=X^1_1X^2_2-X^2_1X^1_2 &\longmapsto & -2\\
X_{4,1}:=(X^1_1X^1_2)^2 &\longmapsto & z^4-1\\
X_{4,2}:=(X^2_1X^2_2)^2 &\longmapsto & \frac{(z^6-z^4+1)^2}{(z^4-1)^3}\\
X_{4,3}:=(X^1_1X^2_2+X^2_1X^1_2)^2 &\longmapsto & \frac{4z^6}{z^4-1} \\
X_{4,4}:=(X^1_1X^1_2)(X^1_1X^2_2+X^2_1X^1_2)&\longmapsto & 2z^3\\
X_{4,5}:=(X^2_1X^2_2)(X^1_1X^2_2+X^2_1X^1_2)&\longmapsto & \frac{2z^3(z^6-z^4+1)}{(z^4-1)^2}\\
X_{4,6}:=X^1_1X^1_2X^2_1X^2_2&\longmapsto & \frac{z^6-z^4+1}{z^4-1},
\end{eqnarray*}
where the arrow refers to the image on the Picard-Vessiot extension under the map
\[
X_{ij}\longmapsto Y^{(i-1)}_j.
\]
We have:
\[
\frac{4X_{4,1}+X_{2,1}^2}{2X_{4,4}}=z.
\]
So the equation is basic. The homogeneous relations that vanishe under this evaluation homomorphism are given by the homogeneous ideal generated by
\begin{eqnarray*}
X_{4,4}X_{4,5}-X_{4,6}X_{4,3}\\
X_{4,6}^2-X_{4,1}X_{4,2}\\
X_{4,1}X_{4,2}-\frac{1}{16}(X_{4,3}-X_{2,1}^2)^2\\
X_{4,3}X_{2,1}^2-(X_{4,3}-2X_{4,6})^2\\
X_{4,4}^2-X_{4,1}X_{4,3}\\
X_{4,5}^2-X_{4,2}X_{4,3}.
\end{eqnarray*}
All of these are invariant under the group $G_F$ generated by $D_\infty$ together with
\[
\left(
\begin{array}{cc}
1 & 0\\
0 & -1
\end{array}
\right),
\]
and so correspond to automorphisms of
the space
\[
\textrm{Proj}(\mathbb{C}[X_{2,1},X_{4,1},X_{4,2},X_{4,3},X_{4,4},X_{4,5},X_{4,6}]).
\]
Again we have an exact sequence:
\[
1\longrightarrow G\longrightarrow G_F\longrightarrow \langle z\mapsto -z \rangle\longrightarrow 1.
\]
\end{exa}

\begin{exa}
We study in more detail the comments in Remark \ref{exa0}. Consider the differential equation ~\cite{ulmer} $L(y)=0$ given by

{\scriptsize
\[
y'''+\frac{3(3x^2-1)}{x(x-1)(x+1)}y''+\frac{221x^4-206x^2+5}{12x^2(x-1)^2(x+1)^2}y'+\frac{374x^6-673x^4+254x^2+5}{54x^3(x-1)^3(x+1)^3}y=0.
\]}
Its Picard-Vessiot extension has differential Galois group $G_{54}$ of order $54$.
The singular points of $L(y)=0$ are $0$, $1$, $-1$ and $\infty$, with respective exponents
\[
\{-\frac{1}{6},\frac{5}{6},\frac{-2}{3}\},\quad\{-\frac{1}{6},\frac{5}{6},\frac{-2}{3}\},
\quad\{-\frac{1}{6},\frac{5}{6},\frac{-2}{3}\},\quad\{\frac{11}{6},\frac{17}{6},\frac{4}{3}\}.
\]
The ramification data in $0$, $1$ and $-1$ is the same, so one can expect some kind of symmetry in between these three points. A quick glance at the equation reveals that all the coefficients of the numerator have even power of $x$, and the denominator present the same exponents for $x-1$ and for $x+1$. This equation admits one symmetry $x\mapsto -x$.

With some computation we can see that if $X$ denotes the solution based at $0$ with exponent $-\frac{1}{6}$, $Y$ the one with $-\frac{5}{6}$ and $Z$ the final one, then
\[
YZ^2+X^3-\frac{16}{81}XY^2=0.
\]
This corresponds to an elliptic curve. The other third degree $G_{54}$-semi-invariant is given by:
\[
XZ^2+\frac{32}{162}X^2Y+\frac{256}{19683}Y^3=\left(\frac{1}{x^3(x^2-1)^3}\right)^{\frac{1}{2}}.
\]
The vanishing $G_{54}$-invariant polynomials are then spanned by
\[
\begin{array}{l}
(YZ^2+X^3-\frac{16}{81}XY^2)(XZ^2+\frac{32}{162}X^2Y+\frac{256}{19683}Y^3)\\
(YZ^2+X^3-\frac{16}{81}XY^2)^2.
\end{array}
\]
Now, as in the case of the equation in the first example over the elliptic curve, we obtain the projection version of $G_{54}$ as a normal subgroup of index two in $F_{36}$. Although it may seem like the theorem applies to this case, the problem here is that the equation is neither basic nor standard. In fact, if we take the symmetric product of this equation with
\[
y'-\frac{2}{3}\left({\frac{1}{x}+\frac{1}{x+1}+\frac{1}{x-1}-\frac{x}{x^2+1}}\right)y=0
\]
we obtain the equation
{\scriptsize
\[
y'''+\frac{5x^4-1}{x(x^2-1)(x^2+1)}y''+\frac{1}{12}\frac{45x^8+20x^6-130x^4+20x^2-3}{x^2(x^2-1)^2(x^2+1)^2}y'
-\frac{20}{27}\frac{x^4-6x^2+1}{x(x^2+1)^3(x^2-1)}y=0.
\]
}
The solutions to this last equation are given by the solutions to the former multiplied by $\sqrt[3]{\frac{(x^3-x)^2}{x^2+1}}$. The Galois group is unmodified, as are the vanishing homogeneous $G_{54}$-invariants. This equation has group of symmetries isomorphic to the dihedral group of order eight:
{\footnotesize
\[
\{x\mapsto x, x\mapsto\frac{1}{x}, x\mapsto-x, x\mapsto\frac{-1}{x}, x\mapsto\frac{-x+1}{x+1}, x\mapsto\frac{x+1}{x-1}, x\mapsto\frac{x+1}{-x+1}, x\mapsto\frac{x-1}{x+1}\},
\]
}
and it descends to the sphere parameterized by $z=\frac{1}{16}\frac{(x^2+1)^4}{x^2(x+1)^2(x-1)^2}$ to the basic equation
\[
y'''+\frac{1}{2}\frac{8z-5}{z(z-1)}y''+\frac{5}{48}\frac{21z-5}{z^2(z-1)}y'-\frac{5}{864}\frac{1}{(z-1)z^3}y=0,
\]
which is expressed in terms of $z$ and $\frac{\partial}{\partial z}$. The quotient group of $F_{36}$ by $PG_{54}$ corresponds to the symmetry coming from lifting this equation to the sphere parameterized by $\sqrt{z}=\frac{1}{4}\frac{(x^2+1)^2}{x(x+1)(x-1)}$, and the remaining symmetries arise from the lifting to the original sphere parameterized by $x$.
\end{exa}

\begin{exa}
Consider the differential equation  ~\cite{singer} $L(y)=0$ given by
\[
y'''+\frac{21(x^2-x+1)}{25x^2(x-1)^2}y'+\frac{21(-2x^3+3x^2-5x+2)}{50x^3(x-1)^3}y=0.
\]
Its Picard-Vessiot extension has differential Galois group $A_5$, the rotational icosahedral group of order 60. This equation admits one symmetry: $z\mapsto -z+1$. The equation is a symmetric power of a second order linear differential equation and so its solutions satisfy the equation
\[
XY-Z^2=0.
\]
The group of automorphisms obtained by conjugating $A_5$ within the subgroup of $GL_3(\mathbb{C})$ fixing the homogeneous ideal generated by this conic is $S_5$ (conjugation with a diagonal matrix with determinant $-1$ together with the inner automorphisms). As in the previous example, because the equation is not standard, the quotient group corresponds to the symmetries of another equation which can be pulled back to obtain this one. Indeed we can take the symmetric product with the equation
\[
y'-\frac{2}{3}\left(\frac{1}{x}+\frac{1}{x-1}\right)y=0
\]
to obtain
{\footnotesize
\[
y'''+\frac{2(2x-1)}{x(x-1)}y''+\frac{1}{75}\frac{163x^2-163x+13}{x^2(x-1)^2}y'-\frac{11}{1350}\frac{2x^3-3x^2-3x+2}{x^3(x-1)^3}y=0,
\]
}
which is an equation with the same Galois group, but with symmetries forming the group $S_3$:
\[
\{x,1-x,\frac{1}{x},\frac{x-1}{x},\frac{1}{1-x},\frac{x}{x-1}\}.
\]
The equation then descends to the sphere parameterized by $z=\frac{4}{27}\frac{(x^2-x+1)^3}{x^2(x-1)^2}$ as the basic equation
\[
y'''+\frac{1}{2}\frac{7z-4}{z(z-1)}y''+\frac{1}{900}\frac{1389z-200}{z^2(z-1)}y'-\frac{11}{5400}\frac{1}{z^2(z-1)}y=0,
\]
which is expressed in terms of $z$ and $\frac{\partial}{\partial z}$. The quotient group of $S_5$ by $A_5$ corresponds to the symmetry arising from lifting this equation to the sphere parameterized by $\frac{(x-c)^3}{x(x-1)}$, where $c=\frac{1+i\sqrt{3}}{2}$; the remaining symmetries come from the lifting to the original sphere parameterized by $x$.
\end{exa}

\section{Comments}

The computations involved in obtaining the Fano group, as well as the normalizer of the differential Galois group, are quite complicated. Among other things, they require extensive use of Van Hoeij and Weil's algorithm ~\cite{We}. In the case where the Galois group is finite, things may be simplified. Indeed, if the charts $(U,z,y^1,\ldots,y^n)$ are taken so that the section $(z,1,0,\ldots,0)$ is cyclic under $\nabla_v$, then the Picard-Vessiot extension is given by $\iota(k)[\frac{\partial\Phi^i}{\partial y^1}]$. So we should be able to replace $R$ with $\mathbb{C}[X^i_1]^G/F_0$, where $F_0$ is the contraction of $F$. This ideal is what Fano originally considered, cf.~\cite{singer88}.


\begin{thebibliography}{99}

\bibitem{baldassarri} F. Baldassarri, On second-order linear differential equations with algebraic solutions on algebraic curves, \emph{Amer. J. Math.} \textbf{102} (1980), no. 3, 517--535.
\bibitem{baldassarri2} F. Baldassarri, B. Dwork, On second order linear differential equations with algebraic solutions, \emph{Amer. J. Math.} \textbf{101} (1979), no. 1, 42--76.
\bibitem{ber} M. Berkenbosch, Algorithms and moduli spaces for differential equations, \emph{S\'eminaires \& Congr\`es}
\textbf{13} (2006) 1–-38.
\bibitem{beukers} F. Beukers, The maximal differential ideal is generated by its invariants, \emph{Indag. Mathem., N.S.} \textbf{11} (1) (2000), 13-18.
\bibitem{churchill} R.C. Churchill, D.L. Rod, B.D. Sleeman, Linear differential equations with symmetries, in Ordinary and partial differential equations Volume V, \emph{Pitman Research Notes in Mathematics}, \textbf{370}, Addison Wesley Longman, Essex, 1997.
\bibitem{compoint} E. Compoint, Differential equations and algebraic relations, \emph{J. symb. Comp.}, \textbf{25} (1998), 705-725.
\bibitem{casale} G. Casale, \textit{Sur le groupo\"ide de Galois d'un feuilletage},
PhD. thesis, Paul Sabatier University 2004, Toulouse, supervised by Prof. d'Emmanuel Paul \& Prof. Jean-Pierre Ramis.
\bibitem{fano} G. Fano, Ueber Lineare Homogene Differentialgleichungen
mit algebraischen Relationen zwischen den Fundamentalloesungen,
\emph{Math. Ann.}, \textbf{53} (1900), 493-590.
\bibitem{vHoeijvdP} M. van Hoeij, M. van der Put, Descent for differential modules and skew fields, \emph{Journal of Algebra} \textbf{296} (2006) 18–-55.
\bibitem{We} M. van Hoeij, J-A. Weil, An algorithm for computing invariants of differential Galois Groups, \emph{J. Pure Appl.
Algebra}, \textbf{117 \& 118} (1997), 353-379.
\bibitem{ince} E.L. Ince, Ordinary Differential Equations, \emph{Dover Publications Inc.} 1958, New York.
\bibitem{magid} A.R. Magid, Lectures on differential Galois theory, \emph{University Lecture Series}, \textbf{7}, American Mathematical Society 1994, Providence, RI.
\bibitem{malgrange} B. Malgrange, Le groupo\"ide de Galois d'un feuilletage, \emph{Monographie} 38
L'enseignement math\'ematique \textbf{2}, 2001.
\bibitem{morales} J.J. Morales Ruiz, Differential Galois theory and non-integrability of Hamiltonian systems,\emph{Progress in Mathematics}, \textbf{179}, Birkhäuser Verlag, Basel, 1999.
\bibitem{nguyen} A.K. Nguyen, A Modern Perspective on Fano's Approach to
Linear Differential Equations, PhD thesis, Groningen, 2008
\bibitem{singer} M. Singer, F. Ulmer, Galois Groups of Second and Third Order Linear Differential Equations, \emph{J. Symb. Comp.}, \textbf{16} (1993), 9-36.
\bibitem{singer88} M.F. Singer, Algebraic relations among solutions
of linear differential equations: Fano's Theorem, \emph{Am. Jour. of
Math.}, \textbf{110} (1988), 115-143.
\bibitem{SvdP} M. van der Put, M.F. Singer, Galois Theory of Linear Differential Equations, \emph{A series of Comprehensive Studies in Mathematics} 328, Springer-Verlag 2003, Berlin Heidelberg New York.
\bibitem{ulmer} F. Ulmer, Liouvillian solutions of third order differential equations, \emph{J. Symb. Comp.}, \textbf{36} (2003), 855-889.
\bibitem{weil} J.A. Weil, First integrals and Darboux polynomials of homogeneous linear differential systems, in Applied algebra, algebraic algorithms and error-correcting codes (Paris, 1995), \emph{Lecture Notes in Comput. Sci.}, \textbf{948}, Springer, Berlin, 1995, 469-484.

\end{thebibliography}
\end{document}